\documentclass[12pt]{article}
\usepackage{amsmath,amssymb,hhline}
\usepackage[russian]{babel}

\usepackage[T2A]{fontenc}    %% for MikTeX  % внутренняя T2A кодировка TeX
\usepackage[cp1251]{inputenc}%% for MikTeX % кодировка - можно использовать [cp866] [cp1251]
\usepackage[russian]{babel} %% for MikTeX % включение переносов

\begin{document}

\begin{center}
{\large\bf R.~M.~Trigub}
\end{center}

\begin{center}
{\large \bf On multiply monotone functions}
\end{center}

\emph{\textbf{Abstract.}} In this paper, the algebra of the differences of two multiply 
monotone functions on $\mathbb{R}_+=(0,+\infty)$ is studied. A sufficient condition for 
the function $f_0\big(|x|_{p,d}\big)$, where
$|x|_{p,d}=\Big(\sum\limits_{j=1}^d|x_j|^p\Big)^{\frac{1}{p}}$,
$p\in(0,+\infty]$, to be represented as the Fourier transform is given. 



\begin{center}
{\large\bf Р.~М.~Тригуб}
\end{center}

\begin{center}
{\large \bf О кратно монотонных функциях}
\end{center}

\emph{\textbf{Аннотация.}} В статье изучена алгебра разностей двух
ограниченных кратно монотонных функций на $\mathbb{R}_+=(0,+\infty)$
и указано достаточное условие представимости функции
$f_0\big(|x|_{p,d}\big)$, где
$|x|_{p,d}=\Big(\sum\limits_{j=1}^d|x_j|^p\Big)^{\frac{1}{p}}$,
$p\in(0,+\infty]$, в виде преобразования Фурье.

Список литературы: 17 названий.

\emph{\textbf{Ключевые слова:}} функции ограниченной вариации,
выпуклые, кратно монотонные, вполне монотонные и положительно
определенные на $\mathbb{R}_+$, преобразование Фурье.

%%% ----------------------------------------------------------------------
\begin{center}
    \textbf{\S0 Введение}
\end{center}

$m$-кратно монотонными называют функции $f:\mathbb{R}_+\rightarrow
\mathbb{R}$, у которых $(-1)^\nu f^{(\nu)}$ $(0\leq \nu\leq m+1)$
принимают значения одного знака (каждая из них и все вместе). Например, убывающая функция с возрастающей
производной (выпуклая)-1-кратно монотонная, а при $m=\infty$ -- вполне монотонная.

 Когда функция $d$ переменных $x_1$, ..., $x_d$ представима на $\mathbb{R}^d$ в
виде преобразования Фурье, т.е. принадлежит винеровской банаховой
алгебре $\Big((x,y)=\sum\limits_{j=1}^d x_j y_j\Big)$
\begin{equation*}
    A(\mathbb{R}^d)=\Big\{f(x)=\widehat{g}(x)=\int\limits_{\mathbb{R}^d}g(y)e^{-i(x,y)}dy,\
    \|f\|_A=\int\limits_{\mathbb{R}^d}|g(y)|dy<\infty\Big\}\quad ?
\end{equation*}

Эта алгебра возникает, напр., при изучении мультипликаторов Фурье из $L_1$ в $L_1$ (см.
[\ref{Stein}--\ref{Stein_Weiss}]).

  Обзор свойств этой алгебры см. в [\ref{Liflyand_Samko_Trigub}]. В
частности, функции из $A(\mathbb{R}^d)$ принадлежат
$C_0(\mathbb{R}^d)$, т.е. непрерывны на $\mathbb{R}^d$ и
\begin{equation*}
    f(\infty)=\underset{\max|x_j|\rightarrow\infty }{\lim f(x)}=0,
\end{equation*}
а также обладают локальным свойством.

Есть много разных достаточных условий принадлежности
$A(\mathbb{R}^d)$ [\ref{Liflyand_Samko_Trigub}]. Особенно хорошо
изучены радиальные функции, т.е. функции вида
$f_0\big(|x|_{2,d}\big)$, где $|x|_{2,d}=\sqrt{(x,x)}$ --- евклидова
норма. В этом случае вопрос о принадлежности $A(\mathbb{R}^d)$ при
$d\geq2$ полностью сводится к принадлежности
$A(\mathbb{R})=A(\mathbb{R}^1)$ другой функции (см.
[\ref{Trigub_Belinsky}], \textbf{6.3.6} и [\ref{Trigub2010}]).

Приведем один пример (см. [\ref{Stein}], гл.~IV, \textbf{7.4}):
\begin{equation*}
    f_0(t)=\frac{e^{it^\alpha}}{(1+t)^\beta},\quad
    f_0\big(|x|_{2,d}\big)\in A(\mathbb{R}^d)\quad
    \Leftrightarrow\quad 2\beta>d\alpha\geq0.
\end{equation*}

Алгебра $A^*(\mathbb{R})$ (ее свойства см. в
[\ref{Belinsky_Lifl_Trig}]) состоит из преобразований Фурье
$f=\widehat{g}$, где $g\in L_1^*(\mathbb{R})$, т.е. удовлетворяет
условию: $\underset{|y|\geq|x|}{{\rm ess}\sup}|g(y)|\in
L_1(\mathbb{R})$.

Beurling [\ref{Beurling}] доказал, что если $f_1(\infty)=0$ и
\begin{equation*}
    \big|f_1(t+h)-f_1(t)\big|\leq\big|f_2(t+h)-f_2(t)\big|\quad (t,h\in\mathbb{R}),
\end{equation*}
где $f_2\in A^*(\mathbb{R})$, то $f_1\in A(\mathbb{R})$.

Отметим еще, что принадлежность $A(\mathbb{R}^d)$ является
существенной при изучении сходимости на $L_1$ и $C$ линейных средних
рядов и интегралов Фурье, определяемых одной функцией--множителем
([\ref{Trigub_Belinsky}], \textbf{8.1.2}), а принадлежность $A^*$
является определяющей для сходимости тех же средних во всех точках
Лебега (почти всюду). См. [\ref{Trigub_Belinsky}], \textbf{8.1.3}.

В настоящей статье изучена алгебра $V_m$ функций, равных разности
двух ограниченных и $m$--кратно монотонных функций на $\mathbb{R}_+$,
и указано достаточное условие для того чтобы
$f_0\big(|x|_{p,d}\big)\in A(\mathbb{R}^d)$, $p\in(0,+\infty]$ (см.
следствия и примеры в \S3).

Далее по плану:

\S1 Кратно монотонные функции.

\S2 Алгебра $V_m(\mathbb{R}_+)$.

\S3 О функциях вида $f_0\big(|x|_{p,d}\big)$ $\big(d\geq2,
p\in(0,+\infty]\big)$.

Через $\gamma(...)$, возможно с индексами, обозначаем положительные
величины, зависящие лишь от переменных, стоящих в скобках.

\begin{center}
    \textbf{\S1 Кратно монотонные функции}
\end{center}

Известно ([\ref{Schoenberg}--\ref{Williamson}]), что если все
функции
\begin{equation*}
    (-1)^\nu f^{(\nu)}\quad (0\leq\nu\leq m-1, m\in\mathbb{N})
\end{equation*}
неотрицательны, убывают и выпуклы вниз на
$\mathbb{R}_+=(0,+\infty)$, то при $t\geq0$ $\big(f(0)=f(+0)\big)$
\begin{equation*}
    f(t)=\int\limits_0^{+\infty}(1-tu)^m_+d\mu(u)\quad
    (a_+=\max\{a,0\}),
\end{equation*}
где $\mu$ --- некоторая положительная мера и конечная на $[0,a]$ при
любом $a>0$.

Известно также, что любая выпуклая (вниз, напр.) функция принадлежит
$AC_{loc}(\mathbb{R}_+)$ (на любом отрезке
$[a,b]\subset\mathbb{R}_+$ абсолютно непрерывна и даже из ${\rm
Lip}1$) и является интегралом от своей (убывающей) правой или левой
производной. Выпуклая функция дифференцируема всюду, кроме,
возможно, не более счетного числа точек, в которых существуют
односторонние производные. Если же выпуклая на $\mathbb{R}_+$
функция ограничена, то она и монотонная.

\textbf{Лемма~1} \emph{Если $f$ ограничена на $\mathbb{R}_+$ и при
$m\geq2$ $\ f^{(m-1)}$ выпукла вверх, то при любом $\nu\in[0,m]$ и
$t>0$ }
\begin{equation*}
    (-1)^{m+\nu+1}f^{(\nu)}(t)\geq0
\end{equation*}
\emph{и существуют пределы $f(+0)$ и $f(+\infty)$. Кроме того,}
\begin{equation*}
    \lim\limits_{t\rightarrow0}t^\nu f^{(\nu)}(t)=\lim\limits_{t\rightarrow+\infty}t^\nu
    f^{(\nu)}(t)=0\qquad (1\leq\nu\leq m)
\end{equation*}
\emph{и }$\displaystyle\int\limits_0^\infty
t^m\big|df^{(m-1)}(t)\big|<\infty$.

\emph{На самом деле, для любой такой ненулевой функции существует
$a\in(0,+\infty]$, при котором и $\nu\geq1$}
\begin{equation*}
    (-1)^{m+\nu+1}f^{(\nu)}(t)>0\qquad \big(t\in(0,a)\big)
\end{equation*}
\emph{и при $a<+\infty$ и $t\geq a$ $\ f(t)=f(+\infty)$. }

$\rhd$ Если функция $f$ ограничена снизу и при некотором
$m\in\mathbb{N}$ \ $f^{(m)}\searrow$ (убывает), то $f^{(m)}(t)\geq0$
при $t>0$. Действительно, как видно из формулы Тейлора, при любом
$a\in \mathbb{R}_+$ и $t\geq a$ (интеграл Стилтьеса)
\begin{equation*}
    f(t)=\sum\limits_{\nu=0}^m\frac{1}{\nu!}f^{(\nu)}(a)(t-a)^{\nu}+\frac{1}{m!}\int\limits_{a}^t(t-a)^mdf^{(m)}(u)\leq\sum\limits_{\nu=0}^m\frac{1}{\nu!}f^{(\nu)}(a)(t-a)^{\nu}.
\end{equation*}

А если бы было $f^{(m)}(a)<0$, то было бы и $f(+\infty)=-\infty$.
Следовательно, существует предел $f^{(m)}(+\infty)$ и
$f^{(m-1)}\nearrow$.

Если $m\geq2$, то по той же причине ($-f^{(m-1)}\searrow$)
$f^{(m-1)}(t)\leq0$ на $\mathbb{R}_+$ и существует предел
$f^{(m-1)}(+\infty)$. При этом $f^{(m)}(+\infty)=0$, так как в
противном случае $f^{(m-1)}$ не может быть ограниченной около
$+\infty$.

Продолжая таким же образом, получаем
\begin{equation*}
    (-1)^{m+\nu+1}f^{(\nu)}(t)\geq0,\qquad f^{(\nu)}(+\infty)=0\quad
    (1\leq\nu \leq m).
\end{equation*}

В силу монотонности функции и ее производных
\begin{equation*}
    \int\limits_0^\infty\sup\limits_{u\geq
    t}|f'(u)|dt=\Big|\int\limits_0^\infty|f'(t)|dt\Big|=\big|f(+\infty)-f(+0)\big|.
\end{equation*}

А если $g(t)\geq0$ и $g\searrow$ на $\mathbb{R}_+$, то при
$\alpha\geq0$
\begin{equation}\label{eq_1}
    0\leq t^{\alpha+1}g(2t)\leq\int\limits_t^{2t}u^\alpha
    g(u)du\rightarrow0\qquad (t\rightarrow+0, t\rightarrow+\infty).
\end{equation}

Поэтому
\begin{equation*}
    \lim\limits_{t\rightarrow+0}tf'(t)=\lim\limits_{t\rightarrow+\infty}tf'(t)=0.
\end{equation*}

Но тогда и
\begin{equation*}
    \int\limits_0^\infty t\sup\limits_{t\geq
    u}|f''(u)|dt=\Big|\int\limits_0^\infty tf''(t)dt\Big|=\Big|\int\limits_0^\infty f'(t)dt\Big|<\infty.
\end{equation*}

В силу \eqref{eq_1}
\begin{equation*}
\lim\limits_{t\rightarrow+0}t^2f''(t)=\lim\limits_{t\rightarrow+\infty}t^2f''(t)=0
\end{equation*}
и т.д.

Еще нужно учесть, что если непрерывная и ненулевая $f^{(\nu)}$ убывает к нулю,
напр., то существует $a_\nu\in(0,+\infty]$, при котором
$f^{(\nu)}(t)>0$ на $(0,a_\nu)$ и $f^{(\nu)}(t)=0$ при $a_\nu<+\infty$ и $t\geq
a_\nu$.

 Очевидно, что $a_{\nu+1}\leq a_{\nu}$. Но
и $a_{\nu}\leq a_{\nu+1}$, так как
\begin{equation*}
    f^{(\nu)}(t)=-\int\limits_t^\infty f^{(\nu+1)}(u)du.
\end{equation*}
\qquad\qquad\qquad\qquad\qquad\qquad
\qquad\qquad\qquad\qquad\qquad\qquad\qquad\qquad\qquad\qquad\qquad\qquad\qquad$\blacktriangleleft$

Вопрос о кратной монотонности стал существенным при определении
положительной определенности, т.е. представлении в виде
преобразования Фурье положительной меры.

 Так, по признаку Пойя, если четная функция $f\in C_0[0,+\infty)$ и
выпукла вниз на $\mathbb{R}_+$, то $f=\widehat{g}$, где $g\in L_1(\mathbb{R})$ и
$g(y)\geq0$.

Более того, $g\in L_1^*(\mathbb{R})$, т.е. $f\in A^*(\mathbb{R})$
(см.[\ref{Trigub_Belinsky}], с.~302).

 А признак типа Пойя для радиальных функций ([\ref{Askey}], см. также
[\ref{Trigub_Belinsky}], \textbf{6.3.7}) теперь можно сформулировать
так:

\emph{если $f_0\in C_0[0,+\infty)$ и при $m=1+\Big[\frac{d}{2}\Big]$
$(d\in\mathbb{N})$\ $(-1)^mf_0^{(m-1)}$ выпукла вверх на
$\mathbb{R}_+$, то}
\begin{equation}\label{eq_2}
    f_0\big(|x|_{2,d}\big)= \int\limits_{\mathbb{R}^d}g(y)e^{-i(x,y)}dy,\quad g\in
    L_1(\mathbb{R}^d),\quad g(y)\geq0\quad (y\in \mathbb{R}^d).
\end{equation}

По теореме Бернштейна ограниченная и вполне монотонная функция на $\mathbb{R}_+$
представима в виде
\begin{equation*}
    f(t)=\int\limits_0^\infty e^{-ut}d\mu(u)\qquad (t\geq0, f(0)=f(+0)),
\end{equation*}
где $\mu$ --- конечная положительная мера $[0,+\infty)$.

А по теореме Schoenberg $f_0\big(|x|_{2,d}\big)$ имеет представление
\eqref{eq_2} при любом $d\in\mathbb{N}$ в том и только в том случае,
когда $f_0(\sqrt{t})$ вполне монотонная.

Отметим еще, что вместе с $m$--кратно монотонной функцией $f$ и суперпозиция $f\circ h$
такая же, если $h(t)\geq0$, $h(t)\not\equiv 0$ и при $t\in \mathbb{R}_+$
\begin{equation*}
    (-1)^{\nu+1}h^{(\nu)}(t)\geq0\qquad (1\leq\nu\leq m).
\end{equation*}

Пример: $h(t)=t^\alpha, \alpha\in(0,1)$.
\newpage
\begin{center}
\textbf{\S2 Алгебра $V_m(\mathbb{R}_+)$}
\end{center}

Разность двух кратно монотонных функций может не быть кратно монотонной.

Обозначим через $V_0(\mathbb{R}_+)$ --- множество функций
ограниченной вариации на $\mathbb{R}_+$, т.е. множество функций,
представимых в виде разности двух ограниченных и монотонных функций.

При $m\in\mathbb{N}$ $V_m(\mathbb{R}_+)$ --- это множество функций с условием
$(f^{(m)}\in V_{0,loc})$
\begin{equation}\label{eq_3}
    \|f\|_{V_m}=\sup\limits_{t\in\mathbb{R}_+}|f(t)|+\int\limits_0^\infty
    t^m\big|df^{(m)}(t)\big|<\infty.
\end{equation}

Это условие при $m\in \mathbb{R}_+$ W.~Trebels [\ref{Trebels73}] использовал как
достаточное условие для мультипликаторов Фурье.Это банахова алгебра. См. также [\ref{Liflyand_Samko_Trigub}].

Функции из $V_1(\mathbb{R}_+)$ называют квазивыпуклыми.

\textbf{Лемма~2} \emph{Для того чтобы $f\in V_1(\mathbb{R}_+)$
необходимо и достаточно, чтобы функция была разностью двух
ограниченных и выпуклых функций на $\mathbb{R}_+$.}

$\rhd$ $\underline{\textrm{Достаточность.}}$

Если $f=f_1-f_2$, где $f_1$ и $f_2$ --- ограниченные выпуклые, то, используя еще
лемму~1,получаем
\begin{equation*}
    \int\limits_0^\infty t|df_{1,2}(t)|=\Big|\int\limits_0^\infty tdf_{1,2}(t)\Big|=\Big|\int\limits_0^\infty f'_{1,2}(t)dt\Big|=\big|f_{1,2}(+0)-f_{1,2}(+\infty)\big|<\infty. \\
\end{equation*}

$\underline{\textrm{Необходимость.}}$

Полагаем $f_1(t)=-\int\limits_0^\infty (u-t)\big|df'(u)\big|\ \
(f_1(+\infty)=0)$.

Очевидно, что при $t>0$
\begin{equation*}
    |f_1(t)|\leq\int\limits_0^\infty u|df'(u)|,\qquad f'_1(t)=\int\limits_t^\infty |df'(u)|\searrow
\end{equation*}
и
\begin{equation*}
    \int\limits_0^\infty t|df'_1(t)|=-\int\limits_0^\infty tdf'_1(t)=\int\limits_0^\infty f'_1(t)dt=-f_1(+0).
\end{equation*}

Так что
\begin{equation*}
    \|f_1\|_{V_1}\leq\int\limits_0^\infty t|df'(t)|+|f_1(+0)|=2\int\limits_0^\infty
    t|df'(t)|.
\end{equation*}

Функция $f_2=f_1-f$ ограничена, как разность ограниченных,
\begin{equation*}
    f'_2(t)=\int\limits_t^\infty|df'(u)|-f'(t)=\int\limits_t^\infty\big(|df'(u)+df'(u)\big)\searrow
\end{equation*}
и
\begin{equation*}
    \|f_2\|_{V_1}\leq \|f_1\|_{V_1}+\|f\|_{V_1}\leq3\|f\|_{V_1}.
\end{equation*}
\qquad\qquad\qquad\qquad\qquad\qquad
\qquad\qquad\qquad\qquad\qquad\qquad\qquad\qquad\qquad\qquad\qquad\qquad$\blacktriangleleft$

Отметим, что функции из $C^2(\mathbb{R}_+)$ образуют плотное
множество в $V_1(\mathbb{R}_+)$. Для доказательства достаточно
применить при $h\rightarrow+0$ функцию Стеклова
\begin{equation*}
    f_{2,h}(t)=\frac {1}{h^2}\int_0^h du_{1} \int _0^h f(t+u_1+u_2)du_2.
\end{equation*}

Введем еще промежуточное пространство между $V_0$ и $V_1$. Это множество
$V_0^*(\mathbb{R}_+)$ функций из $AC_{loc}(\mathbb{R}_+)$  с нормой
\begin{equation*}
    \|f\|_{V_0^*}=\int\limits_0^\infty \underset{u\geq t} {{\rm ess}\sup}|f'(u)|dt.
\end{equation*}

Это банахово пространство не сепарабельное, рефлексивное, в котором
непрерывные функции не образуют плотное множество. Эта алгебра
(кольцо относительно поточечного умножения) существенно отличается
от $L_1(\mathbb{R}_+)$.

Отметим лишь одно отличие $V_0^*$ от $V_1$. Любую функцию из ${\rm Lip}~1$ на отрезке
$[a,b]\subset\mathbb{R}_+$ можно продолжить до функции из $V_0^*(\mathbb{R}_+)$, но не
всегда --- до функции из $V_1$. По этому поводу см. [\ref{Trigub80}], где еще
сравниваются и преобразования Фурье четных функций из $V_0^*$ и $V_1$.

Заметим, что множества $V_1$ и $V_0^*$ можно рассматривать и на отрезке вещественной
оси.

Для отрезка $[0,b]$, напр., появляются нормы
\begin{equation*}
    \int\limits_0^b t\Big(1-\frac{t}{b}\Big)|df'(t)|,\qquad \int\limits_0^b \underset{b\geq u\geq t} {{\rm
    ess}\sup}\big(|f'(u)|+|f'(b-u)|\big)dt.
\end{equation*}

Переходим к общему пространству $V_m$, $m\in\mathbb{N}$.

\textbf{Теорема 1} \emph{Для того чтобы $f\in V_m(\mathbb{R}_+)$
(см. \eqref{eq_3}), необходимо и достаточно, чтобы ее можно было
представить в виде разности двух ограниченных функций с выпуклыми
производными порядка $m-1$.}

$\rhd$ $\underline{\textrm{Достаточность.}}$

Если $f=f_1-f_2$, то используя лемму~1, получаем
\begin{equation*}
    \int\limits_0^\infty t^m\big|df_{1,2}^{(m)}(t)\big|=\Big|\int\limits_0^\infty t^m
    df_{1,2}^{(m)}(t)\Big|=\Big|(-1)^m m!\int\limits_0^\infty
    f'_{1,2}(t)dt\Big|=m!\big|f_{1,2}(+\infty)-f_{1,2}(+0)\big|<\infty.
\end{equation*}

$\underline{\textrm{Необходимость.}}$ Полагаем
\begin{equation*}
    f_1(t)=\frac{(-1)^m}{m!}\int\limits_t^{+\infty}(u-t)^m\big|df^{(m)}(u)\big|.
\end{equation*}

Тогда $f_1$ ограничена:
\begin{equation*}
    0\leq(-1)^m m!f_1(t)\leq\int\limits_0^\infty u^m\big|df^{(m)}(u)\big|
\end{equation*}
и
\begin{equation*}
    f_1^{(m)}(t)=\int\limits_t^\infty\big|df^{(m)}(u)\big|\searrow.
\end{equation*}

При этом
\begin{equation*}
    \|f_1\|_{V_m}\leq\frac{1}{m!}\int\limits_0^\infty u^m\big|df^{(m)}(u)\big|+\int\limits_0^\infty
    t^m\big|df^{(m)}(t)\big|\leq\Big(1+\frac{1}{m!}\Big)\|f\|_{V_m}.
\end{equation*}

И $f_2=f_1-f$ ограничена, как разность двух ограниченных, а
\begin{equation*}
    \begin{split}
        & f^{(m)}_2(t)=\int\limits_t^\infty\big|df^{(m)}(u)\big|-f^{(m)}(t)=\int\limits_t^\infty\big|df^{(m)}(u)\big|+\int\limits_t^\infty df^{(m)}(u)-f^{(m)}(+\infty)\searrow. \\
        &
    \end{split}
\end{equation*}

Предел $f^{(m)}(+\infty)$ существует, так как при $x_1\rightarrow+\infty$ и $x_2>x_1$
\begin{equation*}
    \big|f^{(m)}(x_2)-f^{(m)}(x_1)\big|=\Big|\int\limits_{x_1}^{x_2}
    df^{(m)}(t)\Big|\leq\int\limits_{x_1}^{+\infty}
    t^m\big|df^{(m)}(t)\big|\rightarrow0.
\end{equation*}

При этом
\begin{equation*}
    \|f_2\|_{V_m}\leq\|f_1\|_{V_m}+\|f\|_{V_m}\leq\Big(2+\frac{1}{m!}\Big)\|f\|_{V_m}.
\end{equation*}
\qquad\qquad\qquad\qquad\qquad\qquad
\qquad\qquad\qquad\qquad\qquad\qquad\qquad\qquad\qquad\qquad\qquad\qquad\qquad$\blacktriangleleft$

Trebels [\ref{Trebels75}] (см. также [\ref{Liflyand_Samko_Trigub}], теорема~9.5)
доказал, что при $m>\frac{d-1}{2}$ и $f_0\in V_m(\mathbb{R}_+)$
\begin{equation*}
    f_0(H(x))\in A (\mathbb{R}^d)
\end{equation*}
при любой положительной однородной функции $d$ переменных положительной степени с
условием $H\in C^\infty(\mathbb{R}^d\setminus\{0\})$.

Теперь эту теорему можно сформулировать так:

\emph{если $f_0\in C_0[0,+\infty)$ и ее можно представить в виде разности двух
ограниченных функций, у которых производные порядка $m-1$ при $m>\frac{d-1}{2}$ выпуклы,
то $f_0\circ H\in A(\mathbb{R}^d)$.}
\newpage
\begin{center}
    \textbf{\S3 О функциях вида $f_0\big(|x|_{p,d}\big)$ $(d\geq2, p\in(0,+\infty])$}
\end{center}

При $p\in(0,+\infty)$
\begin{equation*}
    |x|_{p,d}=\Big(\sum\limits_{j=1}^d|x_j|^p\Big)^{\frac{1}{p}},
\end{equation*}
а при $p=\infty$\ $\ \|x\|_{\infty,d}=\max\limits_{1\leq j\leq
d}|x_j|$.

Еще раз об особом случае $p=2$.

\textbf{Предложение} \emph{Если $f_0\in C_0[0,+\infty)$ и $f_0\in
V_m(\mathbb{R}_+)$ при $m=1+\Big[\frac{d}{2}\Big]$, то
$f_0\big(|x|_{2,d}\big)\in A(\mathbb{R}^d)$.}

$\rhd$ Для доказательства достаточно представить $f$ в виде разности согласно теореме~1
и применить признак положительной определенности типа Пойя, приведенный в
\S1.\quad$\blacktriangleleft$

Когда $f_0\big(|x|_{p,d}\big)\in A(\mathbb{R}^d)$ в зависимости от $p$?

Заметим, что при $p\neq2$ не существует частной производной
\begin{equation*}
    \frac{\partial^r f_0\big(|x|_{p,d}\big)}{\partial x_1^r}\qquad (r>p)
\end{equation*}
в точках, в которых $x_2\neq0$. Поэтому лучше учитывать поведение смешанных производных
(дифференцирование по $x_j$ $(j\in[1,d])$ не более одного раза).

\textbf{Теорема~2} \emph{Пусть $f_0\in C_0[0,+\infty)$ и
$f_0^{(d-1)}\in AC_{loc}(\mathbb{R}_+)$. Если еще выполнено условие
А или В, то $f_0\big(|x|_{p,d}\big)\in A(\mathbb{R}^d)$ при
$p\in[1,+\infty]$, а если --- условие Б, то
$f_0\big(|x|_{p,d}\big)\in A(\mathbb{R}^d)$ при $p\in(0,1)$.}

\emph{А. $\int\limits_0^\infty t^{d-1}\underset{u\geq t} {{\rm
    ess}\sup}|f_0(u)|dt<\infty,$}

\emph{Б. $\int\limits_0^\infty t^{dp-1}\underset{u\geq t} {{\rm
    ess}\sup}\ u^{d(1-p)}|f_0^{(d)}(u)|dt<\infty,$}

\emph{В. $f_0(t)=0$ при $t\in[0,a]$, а при $t>0$}
\begin{equation*}
    f_0(t)=O\Big(\frac{1}{t^\varepsilon}\Big),\qquad
    f_0^{(\nu)}(t)=O\Big(\frac{1}{t^{\varepsilon+\nu\delta}}\Big)\quad (\varepsilon>0,
    \delta>1-\frac{2\varepsilon}{d}\Big).
\end{equation*}

\textbf{Следствие 1} \emph{Если $f_0\in C_0[0,+\infty)$ и $f_0\in
V_d(\mathbb{R}_+)$, то при $p\in[1,+\infty]$}

\begin{equation*}
    f_0\big(|x|_{p,d}\big)\in A(\mathbb{R}^d).
\end{equation*}

\textbf{Следствие 2} \emph{Если $f_0\in C_0[0,+\infty)$ и $f_0\in
V_{d+1}(\mathbb{R}_+)$, то при $p\in(0,1)$}

\begin{equation*}
    f_0\big(|x|_{p,d}\big)\in A(\mathbb{R}^d).
\end{equation*}

$\rhd$ Приведем доказательство, пропуская в А и Б некоторые вычисления.

\textbf{Лемма 3} \emph{Если симметричная относительно переменных
$x_j$ $(1\leq j\leq d)$ и чётная по $x_j$ $(1\leq j\leq d)$ функция
$f\in C_0(\mathbb{R}^d)$  имеет на $\mathbb{R}^d_+$ непрерывную
смешанную производную}

\begin{equation*}
    \partial^df(x)=\frac{\partial^d f(x)}{\partial x_1...\partial x_d},
\end{equation*}
\emph{а}
\begin{equation*}
    \lim\limits_{x_1\rightarrow+\infty}\frac{\partial^\nu f(x)}{\partial x_1...\partial
    x_\nu}=0\quad (1\leq\nu\leq d-1),\quad
    \int\limits_{\mathbb{R}^d}\sup\limits_{|x_j|\geq|y_j|,\ 1\leq j\leq
    d}\big|\partial^d
    f(x)\big|dy<\infty,
\end{equation*}
\emph{то $f\in A(\mathbb{R}^d)$.}

Доказательство леммы основано на следующей теореме:

\emph{если для всех $x\in \mathbb{R}^d$ }

\begin{equation*}
    f(x)=\int\limits_{|x_1|}^\infty du_1\int\limits_{|x_2|}^\infty du_2...\int\limits_{|x_d|}^\infty
    g(y_1,...,y_d)dy_d
\end{equation*}
\emph{и}
\begin{equation*}
    \int\limits_{\mathbb{R}^d}\underset{|x_j|\geq|y_j|} {{\rm
    ess}\sup}|g(x)|dy<\infty,
\end{equation*}
\emph{то $f\in A(\mathbb{R}^d)$} ([\ref{Trigub80}], теорема 4,
[\ref{Trigub_Belinsky}], \textbf{6.4.10}).

Эта теорема получена в [\ref{Trigub80}] с использованием полученного там же обобщения
на кратный случай теоремы Бёрлинга, приведенной во введении.

Считая $f_0^{(d)}\in C(\mathbb{R}_+)$, что не уменьшает общности,
полагаем
\begin{equation*}
    g(x)=(-1)^d\partial^d f(x)=(-1)^d\frac{\partial^d f_0\big(|x|_{p,d}\big)}{\partial x_1...\partial
    x_d}.
\end{equation*}

Очевидно, что при $p=\infty$ и $x\in\mathbb{R}^d_+$ $\
\Big|\partial^d
f_0\big(|x|_{\infty,d}\big)\Big|=\Big|f_0^{(d)}\big(|x|_{\infty,d}\big)\Big|$.

Индукцией по $d$ легко доказать, что при $p\in(0,+\infty)$ и
$x\in\mathbb{R}_+^d$
\begin{equation}\label{eq_4}
\frac{\partial^d f_0\big(|x|_{p,d}\big)}{\partial x_1...\partial
    x_d}=\sum\limits_{\nu=1}^d\gamma(d,p,\nu)|x|_{p,d}^{\nu-dp}f_0^{(\nu)}\big(|x|_{p,d}\big)\prod\limits_{j=1}^d
    x_j^{p-1}.
\end{equation}

\textbf{Лемма 4} \emph{Пусть $d\geq2$ и $\alpha>0$. Тогда при
$p=\infty$}

\begin{equation*}
    \int\limits_{\mathbb{R}^d}g\big(|x|_{\infty,d}\big)\prod\limits_{j=1}^d|x_j|^{\alpha-1}dx=\frac{2^d
    d!}{\alpha(\alpha+1)...(\alpha+d-2)}\int\limits_0^\infty
    t^{2\alpha+d-3}g(t)dt,
\end{equation*}
а при $p\in(0,+\infty)$ ($\Gamma$ --- гамма-функция Эйлера).
\begin{equation*}
     \int\limits_{\mathbb{R}^d}g\big(|x|_{p,d}\big)\prod\limits_{j=1}^d|x_j|^{\alpha-1}dx=\frac{2^d
    \Gamma^d\Big(\frac{\alpha}{p}\Big)}{p^{p-1}\Gamma\Big(\frac{d\alpha}{p}\Big)}\int\limits_0^\infty
    t^{d\alpha-1}g(t)dt
\end{equation*}
(в предположении, что простой интеграл справа сходится абсолютно).

При $p\geq1$ $(|x_j|\leq|x|_{p,d}, 1\leq j\leq d)$ и $\nu\geq1$ (см.
\eqref{eq_4})

\begin{equation}\label{eq_5}
    \Big|\frac{\partial^\nu f_0\big(|x|_{p,d}\big)}{\partial x_1...\partial
    x_\nu}\Big|\leq\gamma_0(\nu,p)\max\limits_{1\leq s\leq
    \nu}|x|_{p,d}^{s-\nu}\cdot\Big|f_0^{(s)}\big(|x|_{p,d}\big)\Big|.
\end{equation}

В. Доказательство основано на теореме 2 (Б) из [\ref{Lifl_Trigub}]
(формулируем с учетом симметрии функции):

\emph{если
\begin{equation*}
\frac{\partial^\nu f_0\big(|x|_{p,d}\big)}{\partial x_1...\partial
    x_\nu}=O\Big(\frac{1}{|x|_{2,d}^{\lambda_\nu}}\Big)\quad
    (0\leq\nu\leq d),\ \lambda_0>0
\end{equation*}
и
\begin{equation*}
    \frac{1}{2^d}\sum\limits_{\nu=0}^d\binom{d}{\nu}\lambda_\nu>\frac{d}{2},
\end{equation*}
то $f\in A(\mathbb{R}^d)$.}

В рассматриваемом случае $\lambda_0=\varepsilon$,а при $\nu\geq1$ в силу \eqref{eq_5}
\begin{equation*}
    \lambda_\nu=\nu+\varepsilon+\min\limits_{1\leq
    s\leq\nu}\{\delta-1,s(\delta-1)\}.
\end{equation*}

Так что при $\delta\leq1$ $\ \lambda_\nu=\varepsilon+\delta\nu$, а
при $\delta>1$ $\ \lambda_\nu=\nu+\varepsilon+\delta-1$.

В первом случае $(\delta\leq1)$
\begin{equation*}
    \frac{1}{2^d}\sum\limits_{\nu=0}^d\lambda_\nu\binom{d}{\nu}=\varepsilon+\delta\frac{1}{2^d}\sum\limits_{\nu=1}^d\nu\binom{d}{\nu}=\varepsilon+\delta\frac{1}{2^d}d\cdot2^{d-1}>\frac{d}{2}
\end{equation*}
при $\delta>1-\frac{2\varepsilon}{d}$.

Во втором случае $(\delta>1)$
\begin{equation*}
    \frac{1}{2^d}\sum\limits_{\nu=0}^d\binom{d}{\nu}\lambda_\nu=\frac{1}{2^d}\sum\limits_{\nu=1}^d\binom{d}{\nu}(\nu+\varepsilon+\delta-1)+\frac{\varepsilon}{2^d}=(\varepsilon+(\delta-1))\frac{2^d-1}{2^d}+\frac{d}{2}+\frac{\varepsilon}{2^2}>\frac{d}{2}.
\end{equation*}

Теорема 2
доказана.\qquad\qquad\qquad\qquad\qquad\qquad\qquad\qquad\qquad\qquad\qquad\qquad\qquad$\blacktriangleleft$
\\

$\rhd$ $\underline{\textrm{Доказательство следствия 1.}}$ Если
$f_0\in V_d(\mathbb{R}_+)$, то в силу теоремы 1 она представима в
виде разности $f_1-f_2$ двух ограниченных функций с убывающими
производными порядка $d$. Но тогда (см. еще лемму~1)
\begin{equation*}
    \int\limits_0^\infty t^{d-1}\sup\limits_{t\geq
    u}\Big|f_{1,2}^{(d)}(u)\Big|dt=\int\limits_0^\infty t^{d-1}
    f_{1,2}^{(d)}(t)dt=(-1)^{d-1}(d-1)!\Big(f_{1,2}(+\infty)-f_{1,2}(+ 0\Big).
\end{equation*}

И применяем теорему
2.\qquad\qquad\qquad\qquad\qquad\qquad\qquad\qquad\qquad\qquad\qquad\qquad$\blacktriangleleft$
\\

$\rhd$ $\underline{\textrm{Доказательство следствия 2.}}$ Достаточно
убедиться в неравенстве
\begin{equation*}
    \int\limits_0^\infty t^{dp-1}\sup\limits_{u\geq
    t}u^{d(1-p)}\big|f_0^{(d)}(u)\big|dt\leq\frac{1}{dp}\int\limits_0^\infty
    t^d\big|f_0^{(d+1)}(t)\big|dt.
\end{equation*}

Левая часть не больше
\begin{equation*}
\begin{split}
    &\int\limits_0^\infty t^{dp-1}dt\sup\limits_{u\geq
    t}u^{d(1-p)}\int\limits_u^\infty\big|f_0^{(d+1)}(\nu)\big|d\nu\leq\int\limits_0^\infty
    t^{dp-1}dt\sup\limits_{u\geq t}\int\limits_u^\infty v^{d(1-p)}\big|f_0^{(d+1)}(\nu)\big|d\nu=\\
    &=\int\limits_0^\infty t^{dp-1}dt\int\limits_t^\infty
    u^{d(1-p)}\big|f_0^{(d+1)}(u)\big|du=\int\limits_0^\infty
    u^{d(1-p)}\big|f_0^{(d+1)}(u)\big|du\int\limits_0^u
    t^{dp-1}dt=\\
    &=\frac{1}{dp}\int\limits_0^\infty
    u^d\big|f_0^{(d+1)}(u)\big|du.
\end{split}
\end{equation*}

Осталось повторить доказательство
следствия~1.\qquad\qquad\qquad\qquad\qquad\qquad$\blacktriangleleft$

\textbf{Пример 1}
\begin{equation*}
    f_0(t)=\frac{t^\gamma}{(1+t^\alpha)^\beta},\qquad
    t=|x|_{p,d},\qquad p\in(0,+\infty]
\end{equation*}
принадлежит $A(\mathbb{R}^d)$ только при $\alpha>0$ и
$\alpha\beta>\gamma\geq0$.

$\rhd$ $f_0\in C_0[0,+\infty)$. Поэтому должно быть $\alpha>0$ и
$\alpha\beta>\gamma\geq0$. $f_0\in C^\infty(\mathbb{R}_+)$. Любая
производная $f_0$ в достаточно малой окрестности нуля сохраняет
знак, т.к. при $t\in(0,1)$
\begin{equation*}
    f_0(t)=t^\gamma-\beta t^{\gamma+\alpha}+...
\end{equation*}

Таким же образом ведет себя любая производная и около $\infty$, т.к.
\begin{equation*}
    f_0(t)=t^{\gamma-\alpha\beta}\Big(1-\beta\frac{1}{t^\alpha}+...\Big).
\end{equation*}
\qquad\qquad\qquad\qquad\qquad\qquad\qquad\qquad\qquad\qquad\qquad\qquad\qquad\qquad\qquad\qquad\qquad\qquad$\blacktriangleleft$

\textbf{Пример 2}
\begin{equation*}
    f_0(t)=\frac{e^{it^\alpha}}{(1+t)^\beta},\qquad
    t=|x|_{p,d},\qquad p\in(0,+\infty],\quad \alpha\geq0,\ \beta>0.
\end{equation*}

Если $2\beta>d\alpha$, то $f_0\big(|x|_{p,d}\big)\in
A(\mathbb{R}^d)$ при $p\in [1,\infty]$.А если $\beta>d\alpha$, то $f_0\big(|x|_{p,d}\big)\in A(\mathbb{R}^d )$ при $p\in(0,1)$

$\rhd$ При $p=2$ этот результат точный (см. пример во введении).
Любая производная ${\rm Re} f_0$ и ${\rm Im} f_0$ сохраняет знак в
окрестности нуля и можно, как и в примере~1, применить следствия~1 и
2.

А около $\infty$  применяем при $p\in [1,\infty]$ случай $B$ в теореме 2, учитывая,что
\begin{equation*}
    f_0^{(\nu)}(t)=O\Big(\frac{1}{t^{\beta+\nu(1-\alpha)}}\Big).
\end{equation*}

$\varepsilon=\beta$, $\delta=1-\alpha>1-\frac{2\beta}{d}$ или
$2\beta>d\alpha$.

При $p\in (0,1)$ применяем Б в теореме 2.
\qquad\qquad\qquad\qquad\qquad\qquad\qquad\qquad\qquad\qquad$\blacktriangleleft$

Заметим, что для определения положительной определенности функций
есть результаты, которые существенно зависят от $p\in(0,+\infty]$ (см. напр.,
[\ref{Zastavn}]).

Отметим еще, что, в отличие от случая $d=1$, при $d=2$ средние
арифметические квадратных частных сумм двойного ряда Фурье (суммы
Марцинкевича) могут сходиться не во всех точках Лебега. Это следует
из того, что в $A^*(\mathbb{R}^2)$ нет, практически, функций вида
$f_0\big(|x|_{\infty,2}\big)$ или, что то же самое, функций вида
$f_0\big(|x|_{1,2}\big)$ [\ref{Trigub2015}].

% (см. [\ref{Trigub2000}]).

\end{document}